\documentclass[notitlepage,leqno,11pt]{article}

\usepackage{amscd}
\usepackage{amsmath}
\usepackage{graphicx}
\usepackage{amsfonts}
\usepackage{amssymb}

\newcommand{\proof}[1]{\par\smallskip\noindent{\bf Proof#1.}}
\newcommand{\qed}{\penalty 500\hfill$\square$\par\medskip}

\newcommand{\R}{\mathbb{R}}
\newcommand{\norm}[1]{\left\| #1 \right\|}
\newcommand{\meni}{\leqslant}
\newcommand{\maig}{\geqslant}

\newtheorem{definicion}{Definition}[section]

\newtheorem{theorem}[definicion]{Theorem}

\newtheorem{lemma}[definicion]{Lemma}
\newtheorem{remark}[definicion]{Remark}
\newtheorem{example}[definicion]{Example}
\numberwithin{equation}{section}

\begin{document}

\title{Existence and symmetry for elliptic equations in $\R^n$ with arbitrary growth in the gradient \\ }

\author{ \\ Lucas C. F. Ferreira, Marcelo Montenegro, Matheus C. Santos}
\date{}
\maketitle

\begin{center}
{\small
Universidade Estadual de Campinas, IMECC--Departamento de Matem\'atica, \\ Rua S\'ergio Buarque de Holanda, 651, Campinas-SP, Brazil, CEP 13083-859 \\ \vspace{0.2cm}

E-mail addresses: \texttt{lcff@ime.unicamp.br;msm@ime.unicamp.br;msantos.ufrgs@gmail.com}}

\end{center}

\vspace{0.4cm}

\noindent{\textbf{Abstract}}. We study the semilinear elliptic equation $\Delta u + g(x,u,Du) = 0$ in $\R^n$. The nonlinearities $g$ can have arbitrary growth in $u$ and $Du$, including in particular the exponential behavior. No restriction is imposed on the behavior of $g(x,z,p)$ at infinity except in the variable $x$. We obtain a solution $u$ that is locally unique and inherits many of the symmetry properties of $g$. Positivity and asymptotic behavior of the solution are also addressed. Our results can be extended to other domains like half-space and exterior domains. We give some examples.

\

\noindent{\textbf{AMS Subject Classification 2000}}. 35A01; 35B06; 35B40; 35C15; 35J91.

\

\noindent{\textbf{Keywords}}. Existence; symmetry; asymptotic behavior; semilinear elliptic equations.

\baselineskip 20pt

\section{Introduction}

In this article we study nonlinear elliptic PDEs with the following form
\begin{align}
\Delta u+g(x,u,Du)  &  =0\text{ \ \ \ in \ }\R%
^{n}\label{log1}\\
u  &  \rightarrow0\text{ \ \ as }|x|\rightarrow\infty, \label{log3}%
\end{align}
for $n\geq3$ and $g:\R^n\times\R\times\R^n\rightarrow\R$ verifying $g(x,0,0)\not\equiv 0$ and belonging to a large class of nonlinear functions which include, for example, polynomial and exponential type growths on $u$ or $Du$. Since we are interested in $g$ depending on $u$ and $Du$, we write $g(x,z,p)$ for $z \in \R$, $p \in \R^{n}$ and the gradient of $g$ with respect to the $(n+1)$-last variables will be denoted by $D_{(z,p)}g(x,z,p)$. Throughout the paper, we frequently consider \eqref{log1} with either $g(x,u,Du)$, $g(x,u,|Du|)$, $g(x,|u|,Du)$, or $g(x,|u|,|Du|)$ with the same hypotheses on $g$, except for the symmetry results.

Exponential-type nonlinearities appear naturally in many contexts like conformal geometry and the prescribed curvature problem (see \cite{Cha1}, \cite{Cha2}, \cite{Kaz}), vortex solutions of the Chern-Simons theory (\cite{Caff}, \cite{Chae}, \cite{Str}), statistical mechanics and in a great number of applications as in the description of an isothermal gas sphere and in combustion theory (\cite{G}, \cite{JL}) and stellar structure(\cite{Ch}). On the other hand, nonlinear gradient terms appear naturally in models connected with convective processes and introduces new difficulties when combined with unbounded domains and strong-growth nonlinearities, preventing the use of variational and sub-super solutions methods, Ladyzenskaya-Ural'tseva conditions, Banach fixed point theorem in Sobolev spaces, implicit function theorem, compactness arguments, and Leray-Schauder theory, among others. One of the goals of this work is to provide existence results by using a relatively simpler strategy but new for this prototypical situation.

In smooth bounded domains $\Omega\subset\R^n$, there is a rich literature for \eqref{log1}-\eqref{log3} with general conditions on $g(x,u,Du)$ for existence of solutions, including polynomial or exponential-type growths. In this case existence results have been studied by means of different approaches involving the aforementioned arguments and techniques. For that matter, the reader is referred to \cite{Amann2,Amann3,Cran-Rabi,Choquet-Leray,Kaz2,Kaz3,Quittner1,Rabi1} and their bibliographies. As pointed out in \cite{Choquet-Leray} and \cite{Kaz3}, the use of techniques based on maximum principles in most cases imposes that the nonlinearity grows at most quadratically in $Du$. This kind of restriction appears in the works \cite{Amann2,Bocca-Puel,Dupaigne-Radu,Mignot-Puel}, and was overcomed in \cite{Ruiz} for a logistic equation with $|Du|^q$ with $q>1$ and in bounded domains by combining bifurcation methods and $C^1$-{\it{a priori}} bounds.

For the case of explosive boundary conditions, that is $u\rightarrow\infty\text{ as }x\rightarrow\partial\Omega$ (or as $|x|\rightarrow\infty$), existence of solutions for \eqref{log1} have been addressed in bounded domains $\Omega$ and in $\R^n$ by considering at most polynomial growth at infinity on the gradient $Du$ (see e.g. \cite{Alarcon-Quaas}, \cite{Lions1}, and \cite{Felmer-Quaas-Sirakovc} in $\R^n$). For example, the authors of \cite{Felmer-Quaas-Sirakovc} assumed $-g(u,Du)=f_1(u)\pm f_2(Du)$ with increasing continuous $f$ and $g$ having at most power growth at infinity and $g(x,0,0)=f_1(0)=f_2(0)=0$. We also mention the work \cite{Alarcon-Quaas2} for existence of distributional solutions in $\R^n$ with polynomial growth on both $u$ and $Du$, and without prescribing conditions on $u$ as $|x|\rightarrow\infty$.

Even when $g$ is independent of $Du$, the problem \eqref{log1}-\eqref{log3} in the whole space $\R^n$ with exponential-type growths on $u$ has been considered in dimension $n=2$ in the majority of papers. Usually it is used Trudinger-Moser type inequalities and variational methods for proving existence of solutions (see e.g. results of \cite{Yang} with $n=2$ and its references). In the case of bounded domains, a well known problem arises particularly when
\begin{equation}
g(x,u,Du) = \lambda V(x)e^u, \label{non1}
\end{equation}
which was studied e.g. in \cite{Cran-Rabi,Del,T,Wei-Zhou} (see also their references) with $V$ being a positive bounded smooth function, where the parameter $\lambda$ is assumed to be positive and sufficiently small.

In this work we will show existence of solutions for \eqref{log1}-\eqref{log3} in $\R^n$ with $n \geq 3$ and conditions on $g$ (see \eqref{H1}-\eqref{H3}) covering polynomial and exponential type growths on $u$ and $Du$, see Examples \ref{ex0}, \ref{ex1} and \ref{ex2}. In particular, since $g(x,0,0)$ does not need to be continuous, the nonlinearity \eqref{non1} can be treated with singular potentials $V$ (non-continuous and bounded) and $|\lambda|$ close to zero, including also negative values (see Example \ref{ex1} below). The positivity and symmetry properties of solutions are also addressed, as well as the asymptotic behavior of $u$ and its gradient. By slight modifications on the proofs, our approach can be employed for other unbounded domains like half-space and exterior domains, with either Dirichlet or Neumann homogeneous boundary conditions.

Here we use the integral formulation by means of Green's functions and a contraction argument in the spatial weighted space defined by (for a fixed $k \in \R$)
\begin{equation}
E_{k}\equiv\left\{\;u \mbox{ measurable }\;:\text{ }\text{ess}\sup_{x\in\R^{n}}(1+|x|)^{k}|u(x)|<\infty\right\} \label{espE}%
\end{equation}
and
\begin{equation}
F_{k}\equiv\left\{\;u\in C^1(\R^n)\;:\text{ }\sup_{x\in\R^{n}}
(1+|x|)^{k}\left(|u(x)|+|Du(x)|\right)<\infty\right\}, \label{espD}
\end{equation}
which are Banach spaces with respective norms
\[
\Vert u\Vert_{E_{k}}=\text{ess}\sup_{x\in\mathbb{R}^{n}}(1+|x|)^{k}|u(x)|
\]
and
\[
\Vert u\Vert_{F_{k}}=\sup_{x\in\mathbb{R}^{n}}(1+|x|)^{k}\left(|u(x)|+|Du(x)|\right).
\]
Spaces like above with the homogeneous weight $|x|^{k}$ have been used in \cite{Fer-Mon} to treat the equation $\Delta u + u|u|^{p-2}+V(x)u+f(x)=0$ for $p>n/(n-2)$ with $n \geq 3$. These spaces do not work well for handling nonlinearities with exponential-type growth, mainly because exponential functions transform pole-type singularities into essential ones.

As we will see in the proof of Theorem \ref{t.i}, the choice of a proper value for $k$ in the above spaces depends uniquely on which spaces the function $x\mapsto g(x,0,0)$ is defined and how $|D_{(z,p)}g(\cdot,u,Du)|$ behaves with $|(u,Du)|$.

The problem \eqref{log1}-\eqref{log3} is equivalent to the following
integral equation
\begin{equation}
u(x)=\frac{1}{(n-2)\omega_{n}}\int_{\R^{n}}\frac{1}{|x-y|^{n-2}
}\;g(y,u(y),Du(y))\;dy, \label{int}
\end{equation}
where $\omega_{n}$ is the area of the unit sphere. Therefore, it will be convenient for our purposes to denote the Newtonian potential of a function $f:\R^n\rightarrow \R$ by
\begin{equation}
N(f)(x) := \frac{1}{(n-2)\omega_{n}}\int_{\R^{n}}\frac{1}{|x-y|^{n-2}}\; f(y)\;dy,
\end{equation}
and consider the nonlinear integral operator
\begin{equation}\label{opB}
B(u)(x) := N(g(\cdot,u,Du))(x) =  \frac{1}{(n-2)\omega_{n}}\int_{\R^{n}}\frac{1}{|x-y|^{n-2}}\; g(y,u(y),Du(y))\;dy,
\end{equation}
acting in the space $F_k$.

We shall solve the problem \eqref{int} under the following hypotheses
\begin{eqnarray}
&  g(x,\cdot,\cdot) \mbox{ belongs to } C^1((\R\times\R^n)\backslash(0,0))\cap C(\R\times\R^n), \mbox{ for all } x\in\R^n  \label{H1}; \\
&   \mbox{ There exists } 0<k<n-2 \mbox{ such that the function } \label{H2} \\ & x\mapsto g(x,0,0) \mbox{ belongs to } E_{k+2}; \nonumber\\
&   \mbox{ For the same } k \mbox{ of } (\ref{H2}), \mbox{ there exists } \delta>0 \mbox{ such that } \label{H3} \\ & \displaystyle\sup_{{\tiny 0<\norm{w}_{F_k}\leq\delta}}\norm{D_{(z,p)}g(\cdot,w,Dw)}_{E_2}<\infty,\nonumber
\end{eqnarray}
and a further smallness condition on $sup$ in \eqref{H3}.

\begin{remark}
\label{remark1} For $w\in F_{k}$, $(w,Dw)\equiv 0$ iff $%
w\equiv 0.$ In spite of the fact that $g(x,\cdot ,\cdot )$ is not
differentiable at the point $(0,0)$, we are assuming with \eqref{H3} that $D_{(z,p)}g(x,\cdot ,\cdot )$ is bounded near to the origin. Notice that the supremum of $\left\Vert .\right\Vert _{E_{2}}$ in \eqref{H3} is computed by
excluding $w\equiv 0$.
\end{remark}

The assumptions \eqref{H1}, \eqref{H2} and \eqref{H3} cover many types of nonlinearities with strong growth and gradient dependence. In what follows, we give some examples.
\begin{example} \label{ex0} Recall first that \eqref{log1} is also being defined with $u$ or $Du$ replaced respectively by $|u|$ or $|Du|$ in the arguments of $g$.

\begin{itemize}

\item $g(x,u,|Du|) = \lambda V(x)e^u +\mu W(x)e^{|Du|}$ or $\lambda V(x)e^{e^{{...}^{e^u}}} +\mu W(x)e^{e^{{...}^{e^{|Du|}}}}$, for every $V,W\in E_{k+2}$ with $0<k<n-2$, and $\lambda, \mu \in \R$;

\item $g(x,|u|,|Du|) = W(x)e^{|u|^{m_1}+|Du|^{m_2}}$, $|u|^{m_1}+W(x)e^{|Du|^{m_2}}$, $W(x)e^{|u|^{m_1}}+|Du|^{m_2}$, $W(x)e^{|u|^{m_1}}|Du|^{m_2}+f$, or $W(x)|u|^{m_1} e^{|Du|^{m_2}}+f$, for $m_1,m_2>1$ and $W,f \in E_{k+2}$ with $0<k<n-2$;

\item $g(x,|u|,|Du|) = e^{|u|^{m_1}+|Du|^{m_2}}-1+f(x)$ or $g(x,u,Du) = e^{e^{(|u|^{m_1}+|Du|^{m_2}})}-1+f(x)$, for $m_1,m_2>1$ and $f \in E_{k+2}$ with $0<k<n-2$;

\item $g(x,|u|,|Du|) = |u|^{m_1}+|Du|^{m_2}+f(x)$ or $|u|^{m_1}|Du|^{m_2}+f(x)$, for $m_1,m_2>1$ and $f \in E_{k+2}$ with $0<k<n-2$.
\end{itemize}
\end{example}

Theorem \ref{t.i} corresponds to solving problem \eqref{log1}-\eqref{log3} by looking for a fixed point of the operator $B$ in the space $F_k$ for some suitable choice of $k$, which gives a $C^1$ solution.

A natural question is whether $u$ presents qualitative properties according to $g$. In this direction, if $g$ is symmetric under some orthogonal transformation of $\R^n$, then Theorem \ref{t.iv} guarantees that the solution preserves that symmetry. Also, in Theorem \ref{t.v} we give a condition to improve the natural decay at infinity of the solution belonging to the space $F_k$.

From now on we assume that $n\maig 3$ and that $g:\R^n\times\R\times\R^n\rightarrow\R$ satisfies (\ref{H1})-(\ref{H3}). We begin with existence and local uniqueness of solutions for the integral equation (\ref{int}).

\begin{theorem}
\label{t.i}
There exists a constant $Q_k >0$ such that if $g:\R^n\times\R\times\R^n\rightarrow \R$ satisfies (\ref{H1})-(\ref{H3}) for some $0<k<n-2$ and there is $\varepsilon>0$ such that
\[ G_{\varepsilon}:=\displaystyle\sup_{{\small 0<\norm{w}_{F_k}\leq\varepsilon}}\norm{D_{(z,p)}g(\cdot,w,Dw)}_{E_2} < Q_k
\]
and
\[  \norm{g(\cdot,0,0)}_{E_{k+2}}\meni \varepsilon Q_k,\]
then the integral equation \eqref{int} has a unique solution $u \in F_k$ with $\norm{u}_{F_k}\meni \varepsilon$, which is in particular a weak solution for \eqref{log1}-\eqref{log3}. Furthermore, if $g \in C^{m,\alpha}_{loc}(\R^n\times\R\times\R^n)$ for an integer $m\geq0$ with $0<\alpha<1$, then $u \in C^{m+2,\alpha}_{loc}(\R^n)$ and $u$ verifies \eqref{log1}-\eqref{log3} classically.
\end{theorem}

\begin{remark}
\label{remark10} In the statement of Theorem \ref{t.i}, the constant $Q_k$ can be taken as $\frac{1}{2C_k}$, where $C_k$ is as in Lemma \ref{F} below. See the proof of Theorem \ref{t.i} for more details. In fact, in view of the proof of Lemma \ref{F}, it is possible to estimate $C_k$ and $Q_k$ explicitly.
\end{remark}

In the sequel we present two examples.

\begin{example} \label{ex1}
Let $Q_k=\frac{1}{2C_k}$ where $C_k$ is as in Lemma \ref{F} (see Remark \ref{remark10}). Let $\lambda$ and $\mu$ be real parameters and let
\[g(x,u,Du) = \lambda V(x)e^{u}+\mu W(x)e^{|Du|},\]
where $V,W\in E_{k+2}$ for some $0<k<n-2$. The case $\mu=0$ is the so-called Liouville equation which arises, as pointed out above, in many physical situations and has produced a rich mathematical theory when $n=2$ (see e.g. \cite{bp}, \cite{egp}, \cite{Del}, \cite{T}). Here we solve the problem for all dimension $n \geq 3$. We have that
\begin{eqnarray*}
&(1+|x|)^2|D_{(z,p)} g(x,w,Dw)| = \\
&=  \left( \left(|\lambda|(1+|x|)^2|V(x)|e^{w(x)}\right)^2 + \left(|\mu|(1+|x|)^2|W(x)|e^{|Dw(x)|}\right)^2  \right)^{1/2}  \\
&\meni \left(|\lambda|\norm{V}_{E_{k+2}}+|\mu|\norm{W}_{E_{k+2}}\right) e^{\norm{w}_{F_k}},
\end{eqnarray*}
for all $0\neq w \in F_k$, and
\begin{eqnarray*}
(1+|x|)^{k+2}|g(x,0,0)| &= |\lambda|(1+|x|)^{k+2}|V(x)|+|\mu|(1+|x|)^{k+2}|W(x)| \\
&\meni |\lambda|\norm{V}_{E_{k+2}} + |\mu|\norm{W}_{E_{k+2}}.
\end{eqnarray*}
Then, Theorem \ref{t.i} allows us to solve the problem of the present example if $\lambda$ and $\mu$ satisfy
$$
2C_k\left(|\lambda|\norm{V}_{E_{k+2}}+ |\mu|\norm{W}_{E_{k+2}}\right) e^{2C_k\left(|\lambda|\norm{V}_{E_{k+2}}+ |\mu|\norm{W}_{E_{k+2}}\right)} < 1
$$
and if we take
$$
\varepsilon = 2C_k\left(|\lambda|\norm{V}_{E_{k+2}}+ |\mu|\norm{W}_{E_{k+2}}\right).
$$
The continuous dependence of the solution with respect to $\lambda$ and $\mu$ follows by using that the solution $u$ satisfies
$$
\norm{u}_{F_k}\meni \varepsilon.
$$
with
$$
\varepsilon = 2C_k\left(|\lambda|\norm{V}_{E_{k+2}}+ |\mu|\norm{W}_{E_{k+2}}\right).
$$
This means that the equation $\Delta u + \lambda V(x)e^u + \mu W(x)e^{|Du|} = 0$ has a bounded solution in $\R^n$ if the parameters $|\lambda|$ and $|\mu|$ are small enough, regardless the sign of $\lambda,\mu,V,W$, and allowing to consider non-continuous coefficients $V$ and $W$.
\end{example}

\begin{example}
\label{ex2} According to Remark \ref{remark10}, let us take $Q_k=\frac{1}{2C_k}$ where $C_k$ is as in Lemma \ref{F}. Take $g$ of the form $g(x,z,p_{1},\ldots
,p_{n})=h(x,z^{r_{0}},p_{1}^{r_{1}}\ldots ,p_{n}^{r_{n}})$, where $r_{i}>1$
for all $i$. If $r=\min \{r_{0},\ldots ,r_{n}\}$ and $k=\frac{2}{r-1}$,
suppose that $h(x,0,0)\in E_{k+2}$ and there exists $m>0$ such that $%
D_{(z,p)}h(x,w,Dw)\in E_{m}$ for all $w\in F_{k}$. Then, differentiating we
obtain
\begin{equation*}
D_{(z,p)}g(x,z,p)=\left( r_{0}z^{r_{0}-1}\partial
_{z}h,\;r_{1}p_{1}^{r_{1}-1}\partial _{p_{1}}h,\;\ldots
,\;r_{n}p_{n}^{r_{n}-1}\partial _{p_{n}}h\right) .
\end{equation*}%
If $w\in F_{k}$ with $\left\Vert w\right\Vert _{F_{k}}\leqslant 1$ then
\begin{align*}
& (1+|x|)^{2}|D_{(z,p)}g(x,w,Dw)|\leqslant  \\
& \leqslant \left\vert \left( r_{0}\left[ (1+|x|)^{\frac{2}{r_{0}-1}}|w|%
\right] ^{r_{0}-1}|\partial _{z}h(x,w,Dw)|,\ldots ,r_{n}\left[ (1+|x|)^{%
\frac{2}{r_{n}-1}}|w|\right] ^{r_{n}-1}|\partial _{p_{n}}h(x,w,Dw)|\right)
\right\vert  \\
& \leqslant R\left\vert \left( \left\Vert w\right\Vert _{F_{\frac{2}{r_{0}-1}%
}}^{r_{0}-1}\left\Vert D_{(z,p)}h(\cdot,w,Dw)\right\Vert _{E_{m}},\ldots
,\left\Vert w\right\Vert _{F_{\frac{2}{r_{n}-1}}}^{r_{n}-1}\left\Vert
D_{(z,p)}h(\cdot,w,Dw)\right\Vert _{E_{m}}\right) \right\vert  \\
& \leqslant \sqrt{n+1}R\left\Vert w\right\Vert _{F_{\frac{2}{r-1}}}^{r-1}\left\Vert D_{(z,p)}h(\cdot,w,Dw)\right\Vert _{E_{m}},
\end{align*}
where $R=\max \{r_{0},\ldots ,r_{n}\}$. Thus, for $0<\varepsilon \leqslant 1$,
\begin{equation*}
\displaystyle\sup_{{\small \left\Vert w\right\Vert _{F_{\frac{2}{r-1}%
}}\leq\varepsilon }}\left\Vert D_{(z,p)}g(\cdot ,w,Dw)\right\Vert
_{E_{2}}\leqslant \sqrt{n+1}R\varepsilon ^{r-1}\displaystyle\sup_{{\small
\left\Vert w\right\Vert _{F_{k}}\leq\varepsilon }}\left\Vert D_{(z,p)}h(\cdot
,w,Dw)\right\Vert _{E_{m}}.
\end{equation*}%
If $h$ is such that $\left\Vert h(\cdot,0,0)\right\Vert _{E_{k+2}}\leq \frac{
\varepsilon }{2C_{k}}$ and
\begin{equation*}
\sqrt{n+1}R(2C_{k})^{r}\left\Vert h(\cdot,0,0)\right\Vert _{E_{k+2}}^{r-1}\sup_{{\small \left\Vert w\right\Vert _{F_{k}}\leq\varepsilon }}\left\Vert
D_{(z,p)}h(\cdot ,w,Dw)\right\Vert _{E_{m}}<1,
\end{equation*}
then there exists a solution $u\in F_{k}$ such that $\left\Vert u\right\Vert
_{F_{k}}\leqslant 2C_{k}\left\Vert h(\cdot ,0,0)\right\Vert _{E_{{k+2}}}$.
\end{example}

The solution obtained by the previous theorem inherits many properties from the nonlinearity $g$.

\begin{theorem}\label{t.iii} Under the hypotheses of Theorem \ref{t.i}, the solution $u$ satisfies:
\begin{itemize}
  \item[(i)] If $g \maig 0$, then $u \maig 0$;
  \item[(ii)] If $g(x,z,p)\maig 0$ for all $(x,z,p)\in\R^n\times\R\times\R^n$, with $g(x,z,p)\not\equiv 0$ when $|(z,p)|_{\R\times\R^n}\leq \varepsilon$, then $u>0$;
\item[(iii)] $u$ is radially symmetric provided that $g(\cdot,z,p)$ is radially symmetric for each fixed $(z,p)\in\R\times\R^n$ such that $|(z,p)|_{\R\times\R^n}\leq \varepsilon$.
\end{itemize}
\end{theorem}

More results about symmetry as in item (iii) of Theorem \ref{t.iii} can be proved by considering orthogonal transformations in the space. Let $\mathcal{G}$ be a subset of the orthogonal matrix group $\mathcal{O}(n)$ of $\R^n$. We say that a function $u$ is symmetric under the action of $\mathcal{G}$ when $u(x) = u(Tx)$, for all $T \in \mathcal{G}$. Similarly we say that $u$ is antisymmetric under the action of $\mathcal{G}$ when $u(x) = -u(Tx)$, for all $T \in \mathcal{G}$.

\begin{theorem}\label{t.iv}
Assume the hypotheses of Theorem \ref{t.i} and let $u$ be the solution given by it. Let $\mathcal{G}$ be a subset of $\mathcal{O}(n)$ and suppose that by the action of $\mathcal{G}$, the function $g = g(x,z,p)$ satisfies
\begin{itemize}
  \item[(A)] $g$ is symmetric in $x$ and $p$. Then $u$ is symmetric under $\mathcal{G}$;
  \item[(B)] $g$ is antisymmetric in $p$. Then $u\equiv 0$;
  \item[(C)] $g$ is antisymmetric in $x$, even in $z$ $(i.e. \ g(\cdot,z,\cdot)=g(\cdot,-z,\cdot))$ and symmetric in $p$. Then $u$ is antisymmetric.
\end{itemize}
\end{theorem}

It follows from the definition of the space $F_k$ that the solution given by Theorem \ref{t.i} satisfies $u = \mathcal{O}((1+|x|)^{-k})$ and $Du = \mathcal{O}((1+|x|)^{-k})$ as $|x|\rightarrow \infty$, if $g(x,0,0) = \mathcal{O}((1+|x|)^{-k-2})$. In the next theorem, we improve this behavior by assuming a natural condition, namely if $g(x,0,0) = o((1+|x|)^{-k-2})$ then the solution $u$ and its gradient are $o((1+|x|)^{-k})$ as well.

\begin{theorem}\label{t.v}
Let $g$ be as in Theorem \ref{t.i}. If \ $\lim_{|x|\rightarrow \infty}(1+|x|)^{k+2}|g(x,0,0)| = 0$, then
\begin{equation}\label{limsup}
\lim_{|x|\rightarrow \infty}(1+|x|)^k\left(|u(x)|+|Du(x)|\right) = 0.
\end{equation}
\end{theorem}

In the next section we present the proofs of theorems.

\section{Proof of the Results}

We start by analyzing an integral that will be useful for our needs.
\begin{lemma}\label{lema1}
Let $\alpha,\beta>0$ and $0<n-\alpha<\beta$, then
\[\sup_{x\in\R^n}\int_{\R^n} \frac{1}{|x-y|^\alpha}\frac{1}{(1+|y|)^\beta}\;dy <\infty.\]
\end{lemma}

\proof{}
Using the simplest rearrangement inequality theorem in \cite[p. 82]{Li-Lo}, one has
\[\int_{\R^n} \frac{1}{|x-y|^\alpha}\frac{1}{(1+|y|)^\beta}\;dy \meni \int_{\R^n} \frac{1}{|y|^\alpha}\frac{1}{(1+|y|)^\beta}\;dy\;\;,\;\;\;\forall x \in \R^n ,\]
which is finite, due to the conditions on $\alpha$ and $\beta$.

\qed

The following lemma will be useful for some estimates and its proof can be found in \cite[p. 124]{Li-Lo}.
\begin{lemma}\label{lemaconv}
Let $0<\alpha,\beta<n$ with $0<\alpha+\beta<n$. Then
\[\int_{\R^n}\frac{1}{|y|^{n-\alpha}}\frac{1}{|x-y|^{n-\beta}}\;dy = \frac{C(\alpha,\beta,n)}{|x|^{n-\alpha-\beta}}\]
where $C(\alpha,\beta,n) = \frac{c_\alpha c_\beta c_{n-\alpha-\beta}}{c_{\alpha+\beta} c_{n-\alpha} c_{n-\beta}}$ and $c_\gamma = \pi^{-\gamma/2}\Gamma(\frac{\gamma}{2})$.
\end{lemma}

The next result gives the necessary regularity we will need for the Newtonian potential of a function in the space $E_k$.

\begin{lemma}\label{F}
Let $0<k<n-2$ and $f\in E_{k+2}$. Then $N(f) \in F_k$ and there exists a constant $C_k>0$ satisfying
\begin{equation}\label{C_k}
\norm{N(f)}_{F_k} \meni C_k\norm{f}_{E_{k+2}},\;\;\forall f \in E_{k+2} \ .
\end{equation}
\end{lemma}

\proof{}
First we show that $N(f) \in C^1(\R^n)$. For fixed $x,z \in \R^n$ with $|z|=1$ and $0<t<1/2$, we define the function $h_{y}(s) = |x-y+sz|^{2-n}$ on $[0,t]$. Note that $h_y$ is differentiable on $[0,t]$ if and only if $y \not\in L:=\{x+sz\;|\; s \in [0,t]\}$. If this is the case, we may write
\[
h_y'(s) = \frac{(2-n)z\cdot(x-y+sz)}{|x-y+sz|^{n}}\;\;,\;\; \forall s \in (0,t)
\]
By Mean Value Theorem, for each $y \in \R^n\backslash L$ there exists $t_y \in (0,t)$ such that
\begin{equation}\label{TVM}
\frac{h_y(t)-h_y(0)}{t} = \frac{z\cdot(x-y+t_yz)}{|x-y+t_yz|^{n}}.
\end{equation}
Since $L$ is a measure-zero set, we may write
\begin{eqnarray*}
\frac{N(f)(x+tz)-N(f)(x)}{t} &=& \frac{1}{(n-2)w_n}\int_{\R^n\backslash L} \left(\frac{h_y(t)-h_y(0)}{t} \right)f(y)\;dy \\
&=& -\frac{1}{w_n}\int_{\R^n}\frac{z\cdot(x-y+\overline{t}z)}{|x-y+\overline{t}z|^n}f(y)\;dy.
\end{eqnarray*}
For each $y \in \R^n\backslash L$, let $H_t$ be the function
\[H_t(y) = -\frac{1}{w_n}\frac{z\cdot(x-y+t_yz)}{|x-y+t_yz|^{n}}f(y)\]
where $t_y \in (0,t)$ and satisfies (\ref{TVM}). In spite of the fact that $t_y$ may be not unique, the definition of $H_t(y)$ ensures that a different $t$ satisfying (\ref{TVM}) gives the same value to the expression of $H_t(y)$. Thus $H_t$ is well defined. Furthermore, we have that $H_t\rightarrow H_0$ a.e in $\R^n$.
Note that
\begin{equation}\label{CD1}
|H_t(y)| \meni \frac{1}{w_n}\frac{|f(y)|}{|x-y+t_yz|^{n-1}}\meni G_t(y)
\end{equation}
where
\[G_t(y) = \frac{\norm{f}_{E_{k+2}}}{w_n}\frac{1}{|x-y+t_yz|^{n-1}}\frac{1}{(1+|y|)^{k+2}} \in L^1(\R^n)\;, \mbox{ by Lemma \ref{lema1} }. \]
We also have
\begin{equation}\label{CD2}
G_t(y) \rightarrow G_0(y)\;, \mbox{a.e. in } \R^n \;\mbox { and } \;\int_{\R^n}G_t(y)\;dy = \int_{\R^n}\widetilde{G}_t(y)\;dy
\end{equation}
where
\begin{eqnarray*}
\widetilde{G}_t(y) &=& \frac{\norm{f}_{E_{k+2}}}{w_n}\frac{1}{|y|^{n-1}}\frac{1}{(1+|x+t_yz - y|)^{k+2}}\\
&\meni& \frac{\norm{f}_{E_{k+2}}}{w_n}\frac{1}{|y|^{n-1}}\frac{C_1}{(1+|y|)^{k+2}} \in L^1(\R^n).
\end{eqnarray*}
Therefore, by dominated convergence theorem we have
\[\int_{\R^n}G_t(y)\;dy = \int_{\R^n}\widetilde{G}_t(y)\;dy \rightarrow \int_{\R^n}\widetilde{G}_0(y)\;dy = \int_{\R^n}G_0(y)\;dy.\]
Then, from (\ref{CD1}) and (\ref{CD2}), we conclude that
\[\lim_{t\to 0^+}\frac{N(f)(x+tz)-N(f)(x)}{t} = \lim_{t\to 0^+}\int_{\R^n}H_t(y)\;dy = \int_{\R^n}H_0(y)\;dy.  \]
Thus
\[DN(f)(x)\cdot z = -\frac{1}{w_n}\int_{\R^n}\frac{z\cdot(x-y)}{|x-y|^{n}}f(y)\;dy,\;\;\forall \ |z|=1, \]
and
\[DN(f)(x) = -\frac{1}{w_n}\int_{\R^n}\frac{x-y}{|x-y|^{n}}f(y)\;dy. \]

For a fixed $x_0\in \R^n$ we have
\begin{eqnarray*}
|DN(f)(x_0)-DN(f)(x)| &\meni& \int_{\R^n}\frac{1}{w_n}\left|\frac{x_0-y}{|x_0-y|^{n}}-\frac{x-y}{|x-y|^{n}}\right||f(y)|\;dy
\end{eqnarray*}
and the continuity of $DN(f)$ at $x_0$ follows from the same arguments as above applied to the new functions
\begin{eqnarray*}
H_x(y) &:=& \frac{1}{w_n}\left|\frac{x_0-y}{|x_0-y|^{n}}-\frac{x-y}{|x-y|^{n}}\right||f(y)|;\\
G_x(y) &:=& \frac{\norm{f}_{E_{k+2}}}{w_n}\left(\frac{1}{|x_0-y|^{n-1}}+\frac{1}{|x-y|^{n-1}}\right)\frac{1}{(1+|y|)^{k+2}};\\
\widetilde{G}_x(y) &:=& \frac{\norm{f}_{E_{k+2}}}{w_n}\frac{1}{|y|^{n-1}}\left(\frac{1}{(1+|x_0-y|)^{k+2}}+\frac{1}{(1+|x-y|)^{k+2}}\right)
\end{eqnarray*}
and the estimate
\[ \widetilde{G}_x(y) \meni \frac{C\norm{f}_{E_{k+2}}}{w_n|y|^{n-1}(1+|x_0-y|)^{k+2}} \in L^1(\R^n),\;\mbox{ for }\; |x-x_0|<\frac{1}{2}. \]

For the existence of $C_k$ satisfying (\ref{C_k}), we first note that since
$$
\norm{N(f)}_{F_k} \meni \norm{N(f)}_{E_{k}}+\norm{DN(f)}_{E_{k}}
$$
and the estimates for each term are going to be quite similar, we shall perform only the ones for $DN(f)$.

For $0<k<n-2$, we can apply Lemma \ref{lemaconv} with $\alpha = 1$ and $\beta = n-k-2$ and obtain, for every $x\in\mathbb{R}^n$,
\begin{eqnarray*}
|DN(f)(x)| &\meni& \frac{1}{w_n}\int_{\R^n}\frac{1}{|x-y|^{n-1}}|f(y)|\;dy \\
&=& \frac{1}{w_n}\int_{\R^n}\frac{|y|^{k+2}}{|x-y|^{n-1}}\frac{|f(y)|}{|y|^{k+2}}\;dy \\
&\meni& \frac{1}{w_n}\sup_{y\in\R^n}\left(|y|^{k+2}|f(y)|\right)\int_{\R^n}\frac{1}{|x-y|^{n-1}}\frac{1}{|y|^{k+2}}\;dy \\
&=& \frac{C(n-k-2,1,n)}{w_n}\sup_{y\in\R^n}\left(|y|^{k+2}|f(y)|\right)\frac{1}{|x|^{k+1}} \\
&\meni&\frac{C(n-k-2,1,n)}{w_n}\norm{f}_{E_{k+2}}\frac{1}{|x|^{k+1}}\\
&=:& L_k\norm{f}_{E_{k+2}}\frac{1}{|x|^{k+1}}.
\end{eqnarray*}
Applying Lemma \ref{lema1} with $\alpha = n-1$ and $\beta = k+2$, we conclude
\begin{eqnarray*}
|DN(f)(x)| &\meni& \frac{1}{w_n}\int_{\R^n}\frac{1}{|x-y|^{n-1}}|f(y)|\;dy \\
&=& \frac{1}{w_n}\int_{\R^n}\frac{(1+|y|)^{k+2}}{|x-y|^{n-1}}\frac{|f(y)|}{(1+|y|)^{k+2}}\;dy \\
&\meni& \left(\frac{1}{w_n}\int_{\R^n}\frac{1}{|x-y|^{n-1}}\frac{1}{(1+|y|)^{k+2}}\;dy\right) \norm{f}_{E_{k+2}} \\
&\meni& \left(\frac{1}{w_n}\int_{\R^n}\frac{1}{|y|^{n-1}}\frac{1}{(1+|y|)^{k+2}}\;dy\right) \norm{f}_{E_{k+2}}\\
&=:& M_k\norm{f}_{E_{k+2}}.
\end{eqnarray*}
Therefore, for every $x \in \R^n$,
\begin{eqnarray*}
(1+|x|)^{k+1}|DN(f)(x)| &\meni 2^{k+1}\left(|DN(f)(x)|+|x|^{k+1}|DN(f)(x)|\right)\\
&\meni 2^{k+1}(M_k+L_k)\norm{f}_{E_{k+2}}.
\end{eqnarray*}
Thus $\norm{DN(f)}_{E_{k+1}}\meni 2^{k+1}(M_k+L_k)\norm{f}_{E_{k+2}}$ and by similar calculations we obtain
$$
\norm{N(f)}_{E_{k}}\meni2^{k}(\widetilde{M}_k+\widetilde{L}_k)\norm{f}_{E_{k+2}}
$$
where
\[\widetilde{L}_k = \frac{C(n-2-k,2,n)}{(n-2)w_n} \;\mbox{ and }\; \widetilde{M}_k = \frac{1}{(n-2)w_n}\int_{\R^n} \frac{1}{|y|^{n-2}}\frac{1}{(1+|y|)^{k+2}}\;dy.\]
Thus,
\begin{align*}
\norm{N(f)}_{F_k} &\meni \norm{N(f)}_{E_{k}}+\norm{DN(f)}_{E_{k}} \\
                  &\meni \norm{N(f)}_{E_{k}}+2\norm{DN(f)}_{E_{k+1}} \\
                  &\meni 2^k(M_k + L_k)\norm{f}_{E_{k+2}} + 2^{k+2}(\widetilde{M}_k + \widetilde{L}_k)\norm{f}_{E_{k+2}} \\
                  &\meni 2^{k+2}(\widetilde{M}_k +\widetilde{L}_k+M_k + L_k)\norm{f}_{E_{k+2}}
\end{align*}
and one can take $C_k = 2^{k+2}(M_k + L_k+ \widetilde{M}_k + \widetilde{L}_k)$.

\qed

\proof{ of Theorem \ref{t.i}}
Let $x \in \R^n$, $(x,z_1,p_1),(x,z_2,p_2)\in\R^{n}\times\R\times\R^n$ and $[(z_1,p_1),(z_2,p_2)]$ be the line segment between $(z_1,p_1)$ and $(z_2,p_2)$ in $\R^{n+1}$. If $(0,0)\notin [(z_1,p_1),(z_2,p_2)]$ then, from the hypothesis (\ref{H1}), we have
\[ |g(x,z_1,p_1)-g(x,z_2,p_2)| \meni \sup_{(z,p)\in[(z_1,p_1),(z_2,p_2)]}|D_{(z,p)}g(x,z,p)| |(z_1-z_2,p_1-p_2)| \; .     \]
Now, if $(0,0)\in [(z_1,p_1),(z_2,p_2)]$, then we have that $|(z_1,p_1)|+|(z_2,p_2)| = |(z_1,p_1)-(z_2,p_2)|$ and, by (\ref{H1})
\begin{eqnarray*}
|g(x,z_1,p_1)-g(x,z_2,p_2)| &\meni& |g(x,z_1,p_1)-g(x,0,0)| + |g(x,0,0)-g(x,z_2,p_2)| \\
&\meni& \sup_{(z,p)\in[(z_1,p_1),(z_2,p_2)]\backslash(0,0)}|D_{(z,p)}g(x,z,p)| (|(z_1,p_1)|+|(z_2,p_2)|) \\
&=& \sup_{(z,p)\in[(z_1,p_1),(z_2,p_2)]\backslash(0,0)}|D_{(z,p)}g(x,z,p)| |(z_1-z_2,p_1-p_2)| \; .
\end{eqnarray*}

 Thus, if $u,v \in F_k$, $0< \norm{u}_{F_k},\norm{v}_{F_k}\meni\delta$ and writing $(u,Du)=(u(x),Du(x))$, then
\begin{align*}
|g(x,u,Du)-g(x,v,Dv)| &\meni \sup_{(z,p)\in[(u,Du),(v,Dv)]\backslash(0,0)}|D_{(z,p)}g(x,z,p)| |(u-v,Du-Dv)| \\
&\meni \sup_{{\small 0<\norm{w}_{F_k}\meni\delta}}|D_{(z,p)}g(x,w,Dw)| |(u-v,Du-Dv)|\;.
\end{align*}

Thus,
\[  (1+|x|)^{k+2}|g(x,u,Du)-g(x,v,Dv)| \meni \sup_{{\small 0<\norm{w}_{F_k}\meni\delta}}(1+|x|)^2|D_{(z,p)}g(x,w,Dw)| (1+|x|)^k|(u-v,Du-Dv)|   \]
and by (\ref{H3}), it follows that
\[\norm{g(\cdot,u,Du)-g(\cdot,v,Dv)}_{E_{k+2}} \meni \sup_{{\small0<\norm{w}_{F_k}\meni\delta}}\norm{D_{(z,p)}g(\cdot,w,Dw)}_{E_2}\norm{u-v}_{F_k}   \]

Take $\delta = \varepsilon$ as in the statement of the theorem. We shall show that $B$ is a contraction in the set $A_\varepsilon=\{u \in F_k : \; \norm{u}_{F_k}\meni \varepsilon\}$.

Let $u,v \in A_\varepsilon$ and take $Q_k=\frac{1}{2C_k}$ where $C_k$ is as in Lemma \ref{F}. Noting that $$B(u)-B(v) = N(g(\cdot,u,Du)-g(\cdot,v,Dv)),$$ we can use Lemma \ref{F} and estimate
\begin{align*}
\norm{B(u)-B(v)}_{F_k} &= \norm{N(g(\cdot,u,Du)-g(\cdot,v,Dv))}_{F_k} \\
&\meni C_k\norm{g(\cdot,u,Du)-g(\cdot,v,Dv)}_{E_{k+2}} \\
&\meni C_kG_\varepsilon\norm{u-v}_{F_{k}} \\
&\meni \frac{1}{2}\norm{u-v}_{F_{k}}.
\end{align*}
Thus for $u\in A_\varepsilon$ and $v=0$ in the above inequality, we have
\begin{align*}
\norm{B(u)}_{F_k} &\meni \norm{B(u)-B(0)}_{F_k}+\norm{B(0)}_{F_k} \\
&\meni \frac{1}{2}\norm{u}_{F_{k}} + \norm{N(g(\cdot,0,0))}_{F_k} \\
&\meni \frac{1}{2}\norm{u}_{F_{k}} +C_k\norm{g(\cdot,0,0)}_{E_{k+2}}\\
&\meni \frac{\varepsilon}{2} +C_k\frac{\varepsilon}{2C_k} = \varepsilon
\end{align*}
which shows that $B(A_\varepsilon)\subseteq A_\varepsilon$. Therefore $B$ is a contraction in $A_\varepsilon$ and the result follows by applying the Banach fixed point theorem.

The regularity of $u$ follows from the fact that $u,Du,g(\cdot,u,Du)\in L^\infty(\R^n)$ and the fact that $u$ is a weak solution of \eqref{log1}. Indeed, $u \in W^{1,s}(\Omega)$ and $g(\cdot,u,Du) \in L^s(\Omega)$ for every ball $\Omega$ in $\R^n$ and for every $s>1$, and it solves \eqref{log1} weakly in $\Omega$ without necessarily verifying $u=0$ on $\partial\Omega$. It follows that $u \in W^{2,s}(\Omega)$ for every $s>1$. Then, by the embedding $W^{2,s}(\Omega)\hookrightarrow C^{1,\gamma}(\Omega)$, for $\gamma = 1-\frac{n}{s}$ we conclude that $u \in C^{1,\gamma}(\Omega)$. Therefore, if $g \in C^{m,\alpha}_{loc}(\R^n\times\R\times\R^n)$ then $g(\cdot,u,Du)\in C^{0,\alpha}(\Omega)$ and, by elliptic regularity, we have that $u \in C^{2,\alpha}(\Omega)$. Hence $g(\cdot,u,Du)\in C^{1,\alpha}(\Omega)$ and we can perform the previous argument once more and conclude that $u \in C^{3,\alpha}(\Omega)$. Inductively, we obtain $u \in C^{m+2,\alpha}(\Omega)$, for every ball $\Omega$. In view of the fact that $u$ is a solution of \eqref{log1} in the sense of distributions and $u\in F_k\cap C_{loc}^{m+2,\alpha}(\R^n)$, then $u$ is a classical solution of \eqref{log1}-\eqref{log3}. \qed

\begin{remark}\label{rem-seq1}
The fixed point theorem applied above gives an iterative method to construct the solution $u$, which is the limit in the norm $\norm{.}_{F_k}$ of the following sequence
\[ u_1 = B(0) = N(g(\cdot,0,0)) \;\; \mbox{ and }\;\; u_m = B(u_{m-1})\;, \;\; m \in \mathbb{N}. \]
Moreover, all elements of this sequence verify $\norm{u_m}_{F_k}\leq\varepsilon$.
\end{remark}

\proof{ of Theorem \ref{t.iii}}
The item (i) follows from the fact that the Newtonian potential of a nonnegative function is nonnegative. To prove item (ii), notice that $\norm{u}_{F_k}\leq\varepsilon$ implies that $|(u(x),Du(x))|_{\R\times\R^n}\leq\varepsilon$, for all $x\in\R^n$. It follows that $g(x,u(x),Du(x))\not\equiv 0$, and then $u= N(g(x,u,Du))$ is positive. To establish item (iii), recall first that the solution $u$ is the limit under the norm $\norm{.}_{F_k}$ of the sequence $u_m$ (see Remark \ref{rem-seq1}). Notice that $u_1$ is radially symmetric if and only if $g(x,0,0)$ is radially symmetric. Since $\norm{u_1}_{F_k}\leq\varepsilon$, we have that $|(u_1(x),Du_1(x))|_{\R\times\R^n}\leq\varepsilon$, for all $x\in\R^n$, and then $u_2 = N(g(x,u_1,Du_1))$ is radially symmetric provided that $u_1$ is radially symmetric. By induction, $u_m$ is radially symmetric. Since the convergence in $F_k$ preserves radial symmetry, we conclude that $u$ is radially symmetric.

\qed

\proof{ of Theorem \ref{t.iv}}
(A) Given $T \in \mathcal{G}$, we have that $g(Tx,0,0) = g(x,0,0)$, then
\begin{align*}
u_1(Tx) &= \frac{1}{(n-2)\omega_{n}}\int_{\R^{n}}\frac{1}{|Tx-y|^{n-2}}g(y,0,0)\;dy \\
&= \frac{1}{(n-2)\omega_{n}}\int_{\R^{n}}\frac{1}{|x-T^{-1}y|^{n-2}}g(y,0,0)\;dy \\
&= \frac{1}{(n-2)\omega_{n}}\int_{\R^{n}}\frac{1}{|x-z|^{n-2}}g(Tz,0,0)\;dz \\
&= \frac{1}{(n-2)\omega_{n}}\int_{\R^{n}}\frac{1}{|x-z|^{n-2}}g(z,0,0)\;dz = u_1(x)
\end{align*}
by the change of variables $y = Tz$. Thus, $u_1$ is symmetric under $\mathcal{G}$.

To prove that $u_2$ is symmetric, notice that $Du_1(x) = D (u_1(Tx)) =  T^\top\cdot Du_1(Tx)$. We compute
\begin{align*}
u_2(Tx) &= \frac{1}{(n-2)\omega_{n}}\int_{\R^{n}}\frac{1}{|Tx-y|^{n-2}}g(y,u_1(y),Du_1(y))\;dy \\
&= \frac{1}{(n-2)\omega_{n}}\int_{\R^{n}}\frac{1}{|x-T^{-1}y|^{n-2}}g(y,u_1(y),Du_1(y))\;dy \\
&= \frac{1}{(n-2)\omega_{n}}\int_{\R^{n}}\frac{1}{|x-z|^{n-2}}g(Tz,u_1(Tz),Du_1(Tz))\;dz \\
&= \frac{1}{(n-2)\omega_{n}}\int_{\R^{n}}\frac{1}{|x-z|^{n-2}}g(Tz,u_1(z),T\cdot Du_1(z))\;dz \\
&= \frac{1}{(n-2)\omega_{n}}\int_{\R^{n}}\frac{1}{|x-z|^{n-2}}g(z,u_1(z),Du_1(z))\;dz = u_2(x).
\end{align*}
By the symmetry of $g$. Then $u_2$ is symmetric as well. Using an induction argument, we see that $u_m$ is symmetric under $\mathcal{G}$, for all $m \in \mathbb{N}$. Since $u$ is the limit of $u_m$ in the norm of $F_k$, it preserves the symmetry.
\vspace{0.5cm}

(B) Since $g$ antisymmetric in $p$, then $g(x,0,0) = g(x,0,T0)=-g(x,0,0)$ implies $g(\cdot,0,0)\equiv 0$. Therefore, the fixed point of $B$ is $u \equiv 0$.
\vspace{0.5cm}

(C) One has $g(Tx,0,0) = -g(x,0,0)$, and the computations above give us $u_1(Tx) = -u_1(x)$. Thus, it follows for $u_2$
\begin{align*}
u_2(Tx) &= \frac{1}{(n-2)\omega_{n}}\int_{\R^{n}}\frac{1}{|x-z|^{n-2}}g(Tz,u_1(Tz),Du_1(Tz))\;dz \\
&= \frac{1}{(n-2)\omega_{n}}\int_{\R^{n}}\frac{1}{|x-z|^{n-2}}g(Tz,-u_1(z),T\cdot Du_1(z))\;dz \\
&= -\frac{1}{(n-2)\omega_{n}}\int_{\R^{n}}\frac{1}{|x-z|^{n-2}}g(z,u_1(z),Du_1(z))\;dz = -u_2(x).
\end{align*}
By induction one has $u_m(Tx) = -u_m(x)$. Therefore, one concludes that $u$ is antisymmetric.

\qed

The following lemma is proved in \cite{Fer-Mon}.

\begin{lemma}\label{Lemlimsup}
Let $0<k<n-2$. If $f \in E_{k+2}$, then
\[\limsup_{|x|\rightarrow \infty} |x|^k|N(f)(x)| \meni L_k \limsup_{|x|\rightarrow \infty}|x|^{k+2}|f(x)|\]
where $L_k$ is the constant appearing in the proof of Lemma \ref{F}.
\end{lemma}

\proof{ of Theorem \ref{t.v}}
First recall that the solution given by Theorem \ref{t.i} satisfies $\norm{u}_{F_k}\leq\varepsilon$. Note also that, by the proof of Lemma \ref{F}, if $u \in F_k$, then $g(x,u,Du) \in E_{k+2}$ and therefore $Du = DN(g(x,u,Du))\in E_{k+1}$. Thus, one concludes that
\[\limsup_{|x|\rightarrow \infty}(1+|x|)^{k}|Du(x)| = \limsup_{|x|\rightarrow \infty}\frac{(1+|x|)^{k+1}|Du(x)|}{1+|x|} \meni \limsup_{|x|\rightarrow \infty}\frac{\norm{Du}_{E_{k+1}}}{1+|x|} = 0 \]
Splitting the expression (\ref{limsup}) into two ones, one only needs to check $\lim_{|x|\rightarrow \infty}(1+|x|)^{k}|u(x)|=0$. For that matter, one estimates
\begin{align*}
|g(x,u,Du)| &\meni |g(x,u,Du)-g(x,0,0)| + |g(x,0,0)| \\
&\meni \sup_{{\small 0<\norm{w}_{F_k}\leq\varepsilon}}|D_{(z,p)} g(x,w,Dw)||(u,Du)|+ |g(x,0,0)|.
\end{align*}
Using the hypotheses, one has

\begin{align*}
\limsup_{|x|\rightarrow \infty}|x|^{k+2}|g(x,u,Du)| &\meni \limsup_{|x|\rightarrow \infty}|x|^{k+2}\sup_{{\small0<\norm{w}_{F_k}\leq\varepsilon}}|D_{(z,p)} g(x,w,Dw)||(u,Du)|
\end{align*}
By Lemma \ref{Lemlimsup}, one concludes

\begin{align*}
\limsup_{|x|\rightarrow \infty}|x|^{k}|u(x)| &= \limsup_{|x|\rightarrow \infty}|x|^{k}|B(u)|\\
&= \limsup_{|x|\rightarrow \infty}|x|^{k}|N(g(x,u,Du))| \\
&\meni L_k\limsup_{|x|\rightarrow \infty}|x|^{k+2}|g(x,u,Du)|\\
&\meni L_k\limsup_{|x|\rightarrow \infty}|x|^{k+2}\sup_{{\small0<\norm{w}_{F_k}\leq\varepsilon}}|D_{(z,p)} g(x,w,Dw)||(u,Du)| \\
&\meni L_k\sup_{{\small0<\norm{w}_{F_k}\leq\varepsilon}}\norm{D_{(z,p)} g(\cdot,w,Dw)}_{E_2}\limsup_{|x|\rightarrow \infty}|x|^{k}|(u,Du)| \\
&\meni L_kG_\varepsilon\limsup_{|x|\rightarrow \infty}|x|^{k}(|u(x)| + |Du(x)|) \\
&\meni L_kG_\varepsilon\left(\limsup_{|x|\rightarrow \infty}|x|^{k}|u(x)| + \limsup_{|x|\rightarrow \infty}|x|^{k}|Du(x)|\right)\\
&\meni L_kG_\varepsilon\limsup_{|x|\rightarrow \infty}|x|^{k}|u(x)|
\end{align*}
and, since $L_k G_\varepsilon \meni C_k G_\varepsilon < \frac{1}{2}$, the result follows.

\qed

\noindent{\sc Acknowledgement}. The authors would like to thank an anonymous referee for useful suggestions. L. C. F. Ferreira, M. Montenegro and M. C. Santos have been partially supported by CNPq/Brazil and Capes/Brazil. Part of this work was developed while M. Montenegro was visiting IHP and IH\'ES under CARMIN program and \'Ecole Polytechnique, CMLS, France, their support was greatly appreciated.

\end{document}